\documentclass[10pt]{article}
\usepackage{amssymb}
\usepackage{amsmath, amsthm}
\usepackage{dcolumn}
\usepackage{bm}
\usepackage[T1]{fontenc}
\usepackage{mathptmx}
\usepackage{epsfig} 
\usepackage{graphicx}
\usepackage{epstopdf}

\usepackage[colorlinks=true,citecolor=blue]{hyperref}
\hypersetup{hypertex=true,colorlinks=true,linkcolor=blue,anchorcolor=blue}
\usepackage{cite}
\usepackage{geometry}
\catcode`@=12
\long\def\@makefntext#1{\noindent #1}
\newskip\tabcentering \tabcentering=1000pt plus 1000pt minus 1000pt

\def\MCH#1#2{\setbox0=\hbox{\raise#1\hbox{#2}}\smash{\box0}}

\def\@evenfoot{}\def\@oddfoot{}

\def\sec#1{\vspace{5mm}\leftline{\bf #1}\vspace{3mm}}
\def\subsec#1{\vspace{3mm}\leftline{#1}\vspace{2mm}}

\catcode`@=12

\def\bc{\begin{center}}
	\def\ec{\end{center}}

\def\hang{\hangindent\parindent}
\def\textindent#1{\indent\llap{\qquad #1\ \ \enspace}\ignorespaces}
\def\ref{\par\hang\textindent}

\def\a1{(a_1, a_2, \cdots, a_n)}

\begin{document}
	\thispagestyle{empty}
	\vspace*{-3.0truecm}
	\noindent
	
	\vspace{1 true cm}
	\baselineskip 15pt
	\bc{\large Piston problem to the isentropic  Euler equations for modified Chaplygin gas
		\footnotetext{$^{*}$Corresponding author. \\
		\indent \,\,\,\,\,\,\,\,	E-mail address:  11735032@zju.edu.cn (Meixiang Huang), 454814482@qq.com (Yuanjin Wang), zqshao@fzu.edu.cn (Zhiqiang Shao).\\
	}}\ec

	\vspace*{0.2 true cm}
	\bc{ Meixiang  Huang$^{a, *}$, Yuanjin Wang$^{a}$, Zhiqiang Shao$^{b}$  \\
		{\it $^{a}$School of Mathematics and Statistics,
			Minnan Normal University, Zhangzhou 363000, China	\\
		$^{b}$College of Mathematics and Computer Science, Fuzhou University, Fuzhou 350108, China}
	}\ec

	\vspace*{2.5 true mm}
	\setlength{\unitlength}{1cm}
	\begin{picture}(20,0.1)
	\put(-0.6,0){\line(1,0){14.5}}
	\end{picture}
	
	\vspace*{2.5 true mm}
	\noindent{\small {\small\bf Abstract}
		
\vspace*{2.5 true mm}

In this paper, we solve constructively   the  piston problem for one-dimensional isentropic Euler equations of modified Chaplygin gas. In solutions, we prove rigorously the global existence and uniqueness of a shock wave separating constant states ahead of the piston  when the piston pushed forward into the gas. It is quite different from the results of Chaplygin gas or generalized Chaplygin gas in which a Radon measure solution is constructed to deal with concentration of mass on the piston. When the  piston  pulled back from the gas, we  strictly confirm only the first family  rarefaction wave exists in front of the piston and the concentration will never occur.  In addition, by studying the limiting behavior, we show that the piston solutions of modified Chaplygin gas equations tend to the piston solutions of  generalized or pure Chaplygin gas equations as a single parameter of pressure state function vanishes.

		\vspace*{2.5 true mm}
		
		 \vspace*{2.5 true mm}
		
				\vspace*{2.5 true mm}
		
		\vspace*{2.5 true mm}
		\setlength{\unitlength}{1cm}
		\begin{picture}(20,0.1)
		\put(-0.6,0){\line(1,0){14.5}}
		\end{picture}

\section{Introduction}

~~~~ It is well known that the piston problem is an initial-boundary value problem in fluid dynamics. As shown in [6, 23-24], a thin and long tube with a piston at one end and open at the other end, which is initially filled with static gas, then any movement of the piston will cause the corresponding movement of the gas in the tube. Determining the state of the gas and the propagation of the nonlinear waves in the tube is called a piston problem. 
 In [4, 7] the authors took the one-dimensional piston problem as a model to analyze the occurrence and the motion of the basic nonlinear waves in the isentropic fluid, which shows the importance of the piston problem. In this paper, we  study the one-dimensional piston problem of the isentropic Euler equations [2, 10, 14, 17, 20, 26]

\begin{equation}\label{key}
\left\{\begin{array}{ll} \rho_{t}+(\rho u)_{x}=0,\\
(\rho u)_{t}+(\rho u^{2}+P)_{x}=0,
\end{array}\right . \tag{1.1}
\end{equation} 
where $\rho $,  $u$ and $P$ are the density, velocity and scalar pressure of certain fluid flow located at position $x$ and time $t$, respectively. 
The system (1.1) has been extensively used in numerous fields  due to it's important physical significance. Brenier [2] study the solutions with concentration to the Riemann problem for  the isentropic Chaplygin gas dynamics system. Guo, Sheng and Zhang [11] solved completely this problem, where the delta shock wave solutions were constructed. Roughly speaking, the delta shock wave solution is a solution such that at least one of the variables contains Dirac delta function [3, 5, 21]. In addition, Wang [22] obtained the solutions to the Riemann problem for system (1.1) with the generalized Chaplygin gas, and the formation of delta shock wave was analyzed. Sheng, Wang and Yin [18] clarified the limit behavior of these solutions as the pressure vanishes. Yang and Wang [25] further  studied the system (1.1) for the modified Chaplygin gas, and the limit behavior of constructed Riemann solutions was analyzed as the pressure vanishes.  They also analyzed and identified two kinds of occurrence mechanism on the phenomenon of concentration and the formation of delta shock wave in vanishing pressure limit of solutions to the modified Chaplygin gas equations [26].

As for piston problem, there are numerous excellent papers for related equations and results. Qu, Yuan and Zhao [16] studied the high Mach number limit of the piston problem for the full isentropic Euler equations of polytropic gas. They also  study the one-dimensional piston problem for (1.1) of Chaplygin gas [17], which might be used as an approximate model of polytropic gases in aerodynamics. In their work, a Radon measure solution is constructed when the piston moves forward at sonic or supersonic speed, and   the existence of a shock wave solution is proved as the piston moves forward at subsonic speed, as well as an integral weak solution when the piston recedes from the gas with any constant speed. Fan et.al [10] paid attention to 
the  generalized Chaplygin gas, which substantial difference with Chaplygin gas  lies in that its two characteristic fields are genuinely nonlinear, and
obtained the similar results for the piston problem of the generalized Chaplygin Euler equations. Ding [8] established the global stability of large shock waves to the piston problem for the compressible Magnetohydrodynamics under  small BV perturbations of initial data and the piston speed.

 However, the above works have been focused on the cases that the pressure equation only contains one parameter. It is necessary to investigate the piston problem, in which the pressure equation contains more than one parameter.   
In 2002, the modified Chaplygin gas (MCG) was first proposed by Benaoum [1] to describe the current accelerated expansion of the universe, which is given by

\begin{equation}\label{key}
P = A\rho - \frac{B}{\rho^{\alpha}},  ~~~~~~ 0 < \alpha \leq 1, \tag{1.2}
\end{equation}
where  $A, B > 0$  are two small positive constants, with $A$ representing an ordinary fluid that is subject to linear positive pressure and $B$ relating pressure to some positive power of the inverse for the energy density. It is obvious to see that, when $B = 0$ in (1.2), the pressure is transformed into  the standard equation of state of perfect fluid.  

When $A = 0$, (1.2) reduces into the generalized Chaplygin gas 
\begin{equation}\label{key}
P = - \frac{B}{\rho^{\alpha}},   ~~~~~~ 0 < \alpha \leq 1, \tag{1.3}
\end{equation}
which was extensively studied to compute the lifting force in aerodynamics and analyze dark energy related mathematical problems. Furthermore, as  $\alpha = 1$, (1.3) represents the  Chaplygin gas, 
\begin{equation}\label{key}
P = - \frac{B}{\rho}, \tag{1.4}
\end{equation}
which owns a negative pressure and occurs in certain theories of cosmology. 
Thus, both the Chaplygin and generalized Chaplygin gases are the special cases for the modified Chaplygin gas. Generally, the modified Chaplygin gas can describe the universe to a large extent and explain the mysterious nature of the dark energy and dark matter (see [9, 13]). See the detailed results in  [12, 19] and the references cited therein.

In the present paper, we want to  consider the piston problem for the system (1.1) with modified Chaplygin gas (MCG)
\begin{equation}\label{key}
P = A\rho - \frac{B}{\rho^{\alpha}},  ~~~~~~ 0< \alpha < 1. \tag{1.5}
\end{equation}

In comparison with
 Chaplygin and generalized Chaplygin gases [10, 15, 17], the modified Chaplygin gas has two parameters and can be regard as a combination of the standard fluid and generalized Chaplygin gas, which is more complicated and general; Meanwhile, the analysis of  solvability  for the  piston problem of (1.1) and (1.5) is more complex, especially for the receding case. Recall that the local sound speed of the modified Chaplygin gas is given by
 \begin{equation}\label{key}
c = \sqrt{P'(\rho)} = \sqrt{A+\frac{B \alpha}{\rho^{\alpha+1}}}.  \tag{1.6}
\end{equation}
 Then, the Mach number $M_{0}$ is defined by
 \begin{equation}\label{key}
M_{0} = \frac{ |v_{0}| }{c_{0}},\tag{1.7}
\end{equation}
where $v_{0}$ is  the move speed of piston, $c_{0}$ is as in (1.6) for $\rho = \rho_{0}$. The piston are said to be subsonic, sonic and supersonic when $M_{0} < 1$, $M_{0} = 1$, $M_{0} > 1$, respectively.  In this paper, we discover that there always exists a shock wave  connecting two constant states as the piston rushes into the gas at any constant speed, which quite differs from the results in [10, 17]. Moreover, we prove strictly that only one rarefaction wave is a physically feasible solution, which appears in front of the piston as the piston moves backward from the gas.  Finally, we analyze the limiting behavior of piston solutions to the system (1.1) and (1.5) as the modified Chaplygin gas pressure converges to generalized Chaplygin gas.

The organization of this paper is as follows. In Section 2, we formulate the piston problem of modified Chaplygin gas equations, and shift the coordinates system to move with the piston by Galilean transformation. In Section 3, we briefly recall the results  for the piston problem of (1.1) with (1.3)  as the piston rushes to the gas. In Section 4, the piston solutions to the modified Chaplygin Euler equations of one-dimensional isentropic fluid flow are analyzed and proven in detail as two parameters $A, B>0$, and two  main theorems for the piston problem of (1.1) and (1.5) are presented. In Section 5, when the pressure fades away and $A \rightarrow 0 $, the formation mechanism of singularity is rigorously discussed in the limiting process, and a Radon measure solution is constructed and proven to be a limiting solution of the system (1.1) and (1.5), which is similar with the piston solutions of (1.1) and (1.3).

\section{The piston problem of Modified Chaplygin gas}

~~~~ In this section, we describe the piston problem of modified Chaplygin gas in more details. 
Let the x-axis be the tube, $x = v_{0}t$ is the movement curve of the piston. Suppose the gas fills the domain $\left\lbrace x <0 \right\rbrace $ initially and its velocity and density are $u=0$, $\rho=\rho_{0}$, that is
\begin{equation}\label{key}
U_{0} = (\rho, u) |_{t=0} = (\rho_{0}, 0). \tag{2.1}
\end{equation}
Then the time-space domain is given by
\begin{equation}\label{key}
\Omega_{t} = \left\lbrace (t, x) : x< v_{0}t, ~~~t>0\right\rbrace. \tag{2.2}
\end{equation}
Moreover, the piston is  subject to the usual impermeable condition
\begin{equation}\label{key}
\rho u(t, x) = 0   ~~~~~    on~~~  x = v_{0}t. \tag{2.3}
\end{equation}
The goal of piston problem is to find a solution of (1.1) and (1.5) in the domain $\Omega_{t}$, satisfying (2.1) and (2.3). For the convenience of treating the piston problem,  we use Galilean transformation to shift the coordinates to move with the piston:

\begin{equation}\label{key}
\left\{\begin{array}{ll} 
t' = t,\\
x' = x - v_{0}t,\\
\rho'(t', x') = \rho(t', x'+v_{0}t'),\\
u'(t', x') = u(t', x'+v_{0}t')-v_{0},\\
P'(t', x') = P(t', x'+v_{0}t').
\end{array}\right. \tag{2.4} 
\end{equation}
\noindent
It is easy to verify that the equations in (1.1) are invariant under (2.4), and the domain $\Omega_{t}$ is reduced to $\Omega' = \left\lbrace (t', x') : x'< 0, t'>0\right\rbrace .$ For convenience of statement, we drop all the primes "$'$" and write $\Omega'$ by $\Omega = \left\lbrace (t, x) : x< 0, t>0\right\rbrace $ without confusion. Then the initial condition becomes 
\begin{equation}\label{key}
(\rho, u)|_{t=0} = (\rho_{0}, -v_{0}),  \tag{2.5}
\end{equation}
 and the boundary condition becomes 
\begin{equation}\label{key}
\rho u(t, x) = 0   ~~~~~    on~~~  x = 0. \tag{2.6}
\end{equation}
\noindent
In what follows, we perform the non-dimensional linear transformations of independent and dependent variables, which corresponds to some similarity laws in physics [17]:

\begin{equation}\label{key}
\tilde{t} = \frac{t}{T},~~~~
\tilde{x} = \frac{x}{L},~~~~
\tilde{\rho} = \frac{\rho}{\rho_{0}},~~~~
\tilde{u} = \frac{\sqrt{2}u}{\left| v_{0}\right|},~~~~
\tilde{P} = \frac{2P}{\rho_{0}v_{0}^{2}}, \tag{2.7}
\end{equation}
where $T$ and   $L > 0$ are constants with $\frac{L}{T} = \frac{v_{0}}{\sqrt{2}} $. Inserting (2.7) into (1.1) and simplifying, we can obtain that

\begin{equation}\label{key}
\left\{\begin{array}{ll} \tilde{\rho}_{\tilde{t}} + (\tilde{\rho} \tilde{u})_{\tilde{x}} = 0,\\
 (\tilde{\rho} \tilde{u})_{\tilde{t}} +  (\tilde{\rho} \tilde{u}^{2}+ \tilde{P})_{\tilde{x}} = 0,
\end{array}\right. \tag{2.8} 
\end{equation}
which, together with (2.7), implies that $\tilde{\rho}$, $\tilde{u}$, $\tilde{P}$ still satisfies (1.1) and hence (1.1) is  invariant under (2.7), that is to say, the $\rho$, $u$ and $P$ in (1.1) are equivalent to $\tilde{\rho}$, $\tilde{u}$, $\tilde{P}$ in (2.8). So, from (2.7), in the following we shall take initial data as
\begin{equation}\label{key}
\rho_{0} = 1,~~ v_{0} = \pm \sqrt{2}. \tag{2.9}
\end{equation}
Furthermore, 
by using (1.7) and $c_{0} =  \sqrt{A+\frac{B \alpha}{\rho_{0}^{\alpha+1}}}$, we can obtain
\begin{equation}\label{key}
A \rho_{0}^{\alpha+1} + B \alpha= \frac{\rho_{0}^{\alpha+1} v_{0}^{2}}{M_{0}^{2}}, ~~~~~~
B = \frac{\rho_{0}^{\alpha+1} v_{0}^{2}-A\rho_{0}^{\alpha+1} M_{0}^{2}}{\alpha M_{0}^{2}}.
\tag{2.10}
\end{equation}
Then, (1.5) becomes
\begin{equation}\label{key}
P(\rho) = A\rho - \frac{\rho_{0}^{\alpha+1} v_{0}^{2}-A\rho_{0}^{\alpha+1} M_{0}^{2}}{\alpha M_{0}^{2} \rho^{\alpha}}. \tag{2.11}
\end{equation}
Combining (2.7) with (2.11) and replacing $\frac{\rho_{0}}{\rho}$ by $\frac{1}{\tilde{\rho}}$, we obtain

\begin{equation}\label{key}
\tilde{P}(\tilde{\rho}) = \frac{2A \tilde{\rho}}
{ v_{0}^{2}} - \frac{2(\rho_{0}^{\alpha+1} v_{0}^{2} - A\rho_{0}^{\alpha+1} M_{0}^{2} )}
{\alpha M_{0}^{2} \rho_{0}^{\alpha+1} v_{0}^{2} \tilde{\rho}^{\alpha}}. \tag{2.12}
\end{equation}
Since $\tilde{P}$, $\tilde{\rho}$  in (2.12) are equivalent to corresponding $P$, $\rho$ in (2.11), we conclude that $\rho_{0} = 1$, $ v_{0} = \pm \sqrt{2}$, which   further verifies the consistency of initial data in (2.9). For simplicity of writing, we drop the tildes hereafter. Then,  combining (2.9) and (2.10), one can obtain 

\begin{equation}\label{key}
A  + B\alpha = \frac{2}{M_{0}^{2}}. \tag{2.13}
\end{equation}

In summary, the initial data are defined by
\begin{equation}\label{key}
\rho_{0} = 1,~~ v_{0} = \pm \sqrt{2},~~ P_{0} = A-B =   \frac{2}{ M_{0}^{2}} -B(1+\alpha). \tag{2.14}
\end{equation}

\section{Piston solutions to system (1.1) and (1.3)}
~~~~ For our discussion later, let us briefly recall some results, especially  the construction of Radon measure solutions,  on 
the piston problem of isentropic Euler equations (1.1) with generalized Chaplygin gas (1.3).  
First, we restate the piston problem as 

\begin{equation}\label{key}
	\left\{\begin{array}{ll} \rho_{t}+(\rho u)_{x}=0,\\
		(\rho u)_{t}+(\rho u^{2} - \frac{B}{\rho^{\alpha}})_{x}=0,\\
		(\rho, u)|_{t=0} = (\rho_{0}, -v_{0}),\\
		\rho u(t, x) = 0   ~~~~~    on~~~  x = 0,\\
		\rho_{0} = 1, ~~~~~v_{0} = \pm 1,~~~~~P_{0} = - \frac{1}{\alpha M_{0}^{2}}.\\
	\end{array}\right. \tag{3.1}
\end{equation}
It is easy to know that the piston problem for (1.1) with Chaplygin gas (1.4) is as in (3.1) for $\alpha = 1$.  Notice that the Radon measure solution is not a delta shock, since the singular part is supported on the boundary of the domain.  Now, we recall some definitions about  Radon measure solution [17], which will be employed in this paper.

A Radon measure $m$ on the upper plane $[0, \infty) \times R$ could act on the compactly supported continuous functions
	\begin{equation}\label{key}
	<m, \phi> =  \int_{0}^{\infty}\int_{R}\phi(t,x)m(dxdt),  \tag{3.2}
	\end{equation}
	where the test function $\phi \in C_{0}([0, \infty)\times R)$. The standard  Lebesgue measure $\mathit{L}^{2}$ is one  example of Radon measure on $R^{2}$. The other example is the following Dirac measure supported on a curve [17].

	\noindent \textbf{Definition 3.1.}~~Let $\mathcal{L}$ be a Lipschitz curve given by $x=x(t)$ for  $t$ $  \in$ $ [0,T)$, and $w_\mathcal{L}(t) \in L_{loc}^{1}(0,T)$. The Dirac measure $w_{\mathcal{L}}\delta_{\mathcal{L}}$ supported on $\mathcal{L} \subset R^{2} $ with weight $w_{\mathcal{L}}$ is defined by
	\begin{equation}\label{key}
	<w_{\mathcal{L}}\delta_{\mathcal{L}}, \phi> = \int_{0}^{T}\phi(t, x(t)) w_\mathcal{L}(t) \sqrt{x'(t)^2+1}dt,~~~~~~\forall \phi \in C_{0}(R^2).  \tag{3.3}
	\end{equation}
	
	With this definition, a Radon measure solution can be introduced to construct the solution of (1.1).
	Recall that for two measures $\mu$ and $\nu$, the standard notation $\mu \ll \nu$ means $\nu$ is nonnegative and $\mu$ is absolutely continuous with respect to $\nu$. Now we could formulate the piston problem rigorously by introducing the following definition of measure solutions.

\noindent \textbf{Definition 3.2.}~~ For fixed $0<M_{0}\leq \infty$, let $ \tau, m, n, \wp$ be Radon measures on $\overline{\Omega} = \overline{\{(t, x): x<0,~~ t>0\}}$, and $w_{p}$ a locally integrable nonnegative function on $[0, \infty)$. Then $(\tau, u, w_{p})$ is called a measure solution to the piston problem (3.1) or (3.1) for $\alpha = 1$, provided that

\noindent i) $m\ll \tau$, $n\ll m$, and they have the same Radon-Nikodym derivative $u$; namely
\begin{equation}\label{key}
u \triangleq \dfrac{m(dxdt)}{\tau(dxdt)} = \dfrac{n(dxdt)}{m(dxdt)};
\tag{3.4}
\end{equation}
ii) For any $\phi \in C_{0}^{1}(R^{2})$, there hold
\begin{equation}\label{key}
<\tau, \partial t \phi> + <m, \partial x \phi> + \int_{-\infty}^{0}\rho_{0}\phi(0, x)dx = 0, \tag{3.5}
\end{equation}
\begin{equation}\label{key}
<m, \partial t \phi> + <n, \partial x \phi> +<\wp, \partial x \phi> - <w_{p}\delta_{\{x=0, ~t \geq 0\} }, \phi>
+ \int_{-\infty}^{0}(\rho_{0} u_{0})\phi(0, x)dx = 0, \tag{3.6}
\end{equation}
iii) If $\tau \ll \mathit{L}^2$ with derivative $\rho(t, x)$ in a neighborhood of $(t, x) \in [0, \infty) \times (-\infty, 0]$, and $\wp \ll L^{2}$ with derivative $P(t, x)$ there, then $L^{2}$-a.e. there holds 

\begin{equation}\label{key}
P = -\frac{1}{\alpha \rho^{\alpha}}\frac{1}{M_{0}^{2}} ~~~~~or~~~~~P = -\frac{1}{\rho}\frac{1}{M_{0}^{2}}. \tag{3.7}
\end{equation}
and in addition, the classical entropy condition holds for discontinuities of functions $\rho, u$ near $(t, x)$.

Since the Radon measure solutions exist  only when the piston moves toward the generalized Chaplygin gas (or Chaplygin gas), we just present the corresponding results for the proceeding piston problem in (3.1) and (3.1) for $\alpha = 1$. For the results and relevant proofs of the receding cases, readers can refer to [10, 17], we omit here.

\noindent \textbf{Theorem 3.1.} ~~~ For $v_{0} = -1$, as the piston rushes into the generalized Chaplygin gas (or Chaplygin gas) at Mach number $M_{0} \in (0, \sqrt{\alpha^{-1}})$~ (or $0< M_{0} < 1$),
there exists a shock wave solution of (3.1) (or (3.1) for $\alpha = 1$) connecting states $V_{0} = (1, 1)$ and $V_{1} = (\rho_{1}, 0)$ in the domain $\Omega$, while if Mach number $M_{0} \geq \sqrt{\alpha^{-1}}$ ~(or $M_{0}\geq1$), the Radon measure solution defined above is proved to be a reasonable solution for (3.1) (or (3.1) for $\alpha = 1$) in the domain $\Omega$.

\section{Piston solutions to system (2.8) and (1.5)}

~~~~ In this section, we study the piston problem (2.8), (2.6) and (2.14) for both cases that the piston rushes into and recedes from the MCG (1.5) at a constant speed.  The global existence of solution to piston problem is established and clearly described.

\subsection{Proceeding piston problem} 
~~~~Since a shock wave appears ahead of the piston when it pushes into the gas, we first give the definition of the integral weak solution for problem (2.8), (2.6) and (2.14).

\noindent \textbf{Definition 4.1.}~~ We say $\left(\rho, u \right)\in L^{\infty}(\left[ 0,\infty\right) \times\left( -\infty,0\right] ) $ is an integral solution to  the problem (2.8), (2.6) and (2.14), if for any $\phi \in C^{1}_{0}(R^{2})$, there hold

\begin{equation}\label{key}
\left\{\begin{array}{ll} \int_{\Omega}(\rho \partial t \phi+\rho u \partial x \phi)dxdt + \int_{-\infty}^{0}\rho_{0}(x)\phi(0,x)dx = 0,\\
\int_{\Omega}(\rho u \partial t \phi ++( \rho u^{2}+P) \partial x \phi)dxdt -
\int_{0}^{\infty}P(t,0) \phi(t,0)dt 
+ \int_{-\infty}^{0}\rho_{0}(x)u_{0}(x) \phi(0,x)dx = 0. \tag{4.1}
\end{array}\right .
\end{equation}

Noticing that the problem (2.8), (2.6) and (2.14) is a Riemann problem with boundary conditions for fixed $M_{0} \in (0, \infty)$. Since the equations (2.8), the initial data (2.14) as well as the boundary condition (2.6) are invariant by the transformation 
$$
(x, t)\rightarrow (\mu x, \mu t)   ~~~~~for ~~~\mu  \neq 0.
$$
Thus, we can construct a piecewise constant
self-similar solution $U(x, t) = V(\frac{x}{t}) $ to connect two states $(\rho_{0}, -v_{0})$ and $(\rho_{1}, 0)$, which is in the form

\begin{equation}\label{key}
	U(x, t) = V\left( \frac{x}{t}\right) =
	\left\{\begin{array}{ll} 
		V_{0} = (1, \sqrt{2}),~~~~~ -\infty \leq \frac{x}{t} < \sigma,\\
		V_{1} = (\rho_{1}, 0),~~~~~  \sigma < \frac{x}{t} \leq 0, \tag{4.2}
	\end{array}\right .
\end{equation}
where $V_{0}$ and $V_{1}$ are subject to the Rankine-Hugoniot condition

\begin{equation}\label{key}
\left\{\begin{array}{ll} 
\sigma(\rho_{1} - \rho_{0}) = \rho_{1}u_{1} - \rho_{0}u_{0},~~~~\\
\sigma(\rho_{1}u_{1} - \rho_{0}u_{0}) = \rho_{1}u_{1}^{2} + P_{1} - \rho_{0}u_{0}^{2} - P_{0} .
\end{array}\right. \tag{4.3}
\end{equation}
\noindent
In view of $\rho_{0} = 1$, $u_{0} = \sqrt{2}$, $u_{1} = 0$, it follows from $(4.3)_{1}$ that

\begin{equation}\label{key}
\sigma = -\frac{\sqrt{2}}{\rho_{1}-1}. \tag{4.4}
\end{equation}
\noindent
Note that $\sigma <0$ requires that $\rho_{1} >1$. Substituting $\sigma$ into $(4.3)_{2}$ gives
\begin{equation}\label{key}
P_{1} = \frac{2}{\rho_{1}-1} + 2 + P_{0} 
= \frac{2}{\rho_{1}-1}+2+\frac{2}{M_{0}^2}-B(1+\alpha). \tag{4.5}
\end{equation}
\noindent
Replacing $P_{1}$ in (4.5) by $P_{1} = \left( \frac{2}{M_{0}^2}-B\alpha \right) \rho_{1}
-\frac{B}{\rho_{1}^\alpha}$, then we obtain

\begin{equation}\label{key}
\frac{2(\rho_{1}-1)}{M_{0}^2}=\frac{2\rho_{1}}{\rho_{1}-1}
+B\alpha(\rho_{1}-1)+B\left( \frac{1}{\rho_{1}^\alpha}-1\right) . \tag{4.6}
\end{equation}
From (2.13), we infer that $A$ and $B$ have the form of $\frac{x}{M_{0}^2}$, in which the $x$ is unknown. Without loss of generality, we assume $x$ is a constant independent of $M_{0}$ as
$A$ and $B$ are constants. Then, we know $AM_{0}^2$ and $BM_{0}^2$
are constants independent of $M_{0}$. Multiplying both sides of the equation (4.6) by $M_{0}^2$, we have
\begin{equation}\label{key}
2\left( \rho_{1}-1\right)  = \frac{2\rho_{1}M_{0}^2}{\rho_{1}-1}
+BM_{0}^2\alpha(\rho_{1}-1)+BM_{0}^2\left(\frac{1}{\rho_{1}^\alpha}-1 \right) . \tag{4.7}
\end{equation}
By a simple calculation, (4.7) is equivalent to

\begin{equation}\label{key}
M_{0}^2 = \frac{\left[ (2-BM_{0}^2\alpha)(\rho_{1}-1)-BM_{0}^2\left(\frac{1}{\rho_{1}^\alpha}-1 \right) \right] (\rho_{1}-1)
}{2\rho_{1}}= \frac{\left[ AM_{0}^2(\rho_{1}-1)-BM_{0}^2\left(\frac{1}{\rho_{1}^\alpha}-1 \right) \right] (\rho_{1}-1)
}{2\rho_{1}}, \tag{4.8}
\end{equation}
In what follows, we define a continuous function $f(\rho) = \frac{\left[ AM_{0}^2(\rho-1)-BM_{0}^2\left(\frac{1}{\rho^\alpha}-1 \right) \right] (\rho-1)
}{2\rho}$ with respect to $\rho$. Then,  we can get $f(1)=0$ and $f(\infty)=\infty$.  According to the intermediate value theorem of continuous function,  there exists a $\rho_{1} > 1$, such that $0 < f(\rho_{1}) < \infty$.
Furthermore, we can calculate
\begin{equation}\label{key}
\begin{split}
f'(\rho) &=  \frac{1}{2}\left\lbrace \left( AM_{0}^2+\alpha BM_{0}^2 \frac{1}{\rho^{\alpha+1}}	\right) \left(1-\frac{1}{\rho} \right)+
\left[ AM_{0}^2(\rho-1)-BM_{0}^2\left( \frac{1}{\rho^{\alpha}}-1\right) 
\right]  \frac{1}{\rho^{2}}		 \right\rbrace  \\&=\frac{1}{2} \left[  AM_{0}^2 \left(1-\frac{1}{\rho^{2}}\right)  +\alpha BM_{0}^2  \left( \frac{1}{\rho^{\alpha+1}}-\frac{1}{\rho^{\alpha+2}}
\right)+BM_{0}^2
\left( \frac{1}{\rho^2}-\frac{1}{\rho^{\alpha+2}}
\right) 
\right] .
\end{split}\tag{4.9}
\end{equation}
From (4.9), it is easy to infer $f'(\rho)>0$ for $\rho>1$, i.e., $f(\rho)$ is always a monotone and increasing function for $\rho>1$, 
which guarantees the uniqueness of the shock wave solution. Thus, we  solve  constructively the existence and uniqueness of the shock wave solution for  the proceeding piston problem (2.8), (2.6) and (2.14). In other word, if (2.13) holds, for any $A$ and $B$, 
there alway have a shock wave solution connecting states $V_{0} $ and $V_{1}$ for $0<M_{0} <\infty$, which is quite different from the results in [10, 17].

 Therefore, we have the following result.

\noindent \textbf{Theorem 4.1.} ~~ For $v_{0} = -\sqrt{2}$ and any fixed $A$ and $B$, as the piston rushes into the MCG at arbitrary Mach number $M_{0} \in (0, \infty)$, there always exists a unique shock wave solution of the system (2.8), (2.6) and (2.14) connecting states $V_{0} = (1, \sqrt{2})$ and $V_{1} = (\rho_{1}, 0)$ in the domain $\Omega$ due to the first term of (1.5), which is quite different from the results of Chaplygin gas and generalized Chalygin gas. 

 \noindent \textbf{Remark 4.1.}~~~For any $A$ and $B$, we have demonstrated the global existence and uniqueness of shock wave solution in the one-dimensional piston problem for (2.8), (2.6) and (2.14) as the piston rushes into the MCG for $0<M_{0} <\infty$. We discover that,   the results of (2.8), (2.6) and (2.14) are mathematically different to the isentropic Euler equations with pure or generalized Chaplygin gas, in which a Radon measure solution shall be considered to deal with concentration of mass on the piston   when the Mach numbers satisfy $M_{0}\geq1 $ or $M_{0}\geq \sqrt{\alpha^{-1}}$. In essence, the difference is caused by the first term of (1.5). In  Section 5, we will analyze the formation of singularity on the surface of the piston as $A \rightarrow 0$, which shows that  a Radon measure solution is the limiting piston solution to the modified Chaplygin gas equations as $A \rightarrow 0$.

\subsection{Receding piston problem}

~~~~Let $U(x, t) = (\rho(x, t), u(x, t))$. The system (1.1) with (1.5) have two eigenvalues

 \begin{equation}\label{key}
 \lambda_{1}(U) = u - \sqrt{A + \frac{B \alpha}{\rho^{\alpha+1}}}, ~~~~~ \lambda_{2}(U) = u +  \sqrt{A + \frac{B \alpha}{\rho^{\alpha+1}}}, \tag{4.10}
 \end{equation}
 with the corresponding right eigenvectors
 
 \begin{equation}\label{key}
 \vec{r}_{1} = \left(  \sqrt{A + \frac{B \alpha}{\rho^{\alpha+1}}}, -\rho\right)^{T},~~~~~
  \vec{r}_{2} = \left(  \sqrt{A + \frac{B \alpha}{\rho^{\alpha+1}}}, \rho\right)^{T}, \tag{4.11}
 \end{equation}
satisfying 
\begin{equation}\label{key}
\triangledown \lambda_{i}(U) \cdot\vec{r}_{i} = \frac{2A\rho^{\alpha}+B\alpha(1-\alpha)}{2\rho^{\alpha+1} \sqrt{A + \frac{B \alpha}{\rho^{\alpha+1}}}} \neq 0 ~~~ (i=1, 2). \tag{4.12}
\end{equation}
Thus, we can see that system  (1.1) with (1.5)  is strictly hyperbolic and both characteristic fields are genuinely nonlinear when A, B > 0, in which the associated waves are either rarefaction waves or shock waves [26]. 

\noindent
Now, we can as well construct a solution of the form $U(t,x) = V\left( \frac{x}{t}\right) $. For any fixed $0<M_{0}<\infty$, we suppose the solution is composed of two constant states $V_{0} = (1, -\sqrt{2})$, $V_{1} = (\rho_{1}, 0)$, and a rarefaction wave $V_{m}$ connecting them:

\begin{equation}\label{key}
	U(x, t) = V\left( \frac{x}{t}\right)  =
\left\{\begin{array}{ll} 
V_{0} ,~~~~~~~~~~~~~~~ -\infty \leq \frac{x}{t} < \lambda_{i}(V_{0}),\\

V_{m}\left(\frac{x}{t} \right)  ,~~~~~~~~~~ \lambda_{i}(V_{0}) \leq \frac{x}{t} < \lambda_{i}(V_{1}),\\

V_{1},~~~~~~~~~~~~~~~~~~  \lambda_{i}(V_{1}) < \frac{x}{t} \leq 0,
\end{array}\right. ~~~ (i=1, 2).  \tag{4.13}
\end{equation}

\noindent If the piston recedes from the gas, then $v_{0} >0$ and $\rho_{1} < \rho_{0}$. For the first family rarefaction wave $R_{1}(U)$, we can deduce the self-similar solution  as follow: 

\begin{equation}\label{key}
\left\{\begin{array}{ll} 
\eta = \frac{x}{t} = \lambda_{1}(U) = u -   \sqrt{A + \frac{B \alpha}{\rho^{\alpha+1}}},\\

u -\frac{2}{\alpha+1}
\sqrt{A+\frac{B\alpha}{\rho^{\alpha+1}}} + \frac{\sqrt{A}}{\alpha+1}
\ln\left( 2\sqrt{A}\rho^{\alpha+1}
\left( \sqrt{A+\frac{B\alpha}{\rho ^{\alpha+1}}} + \sqrt{A}\right) +B\alpha \right)  =\\
u_{0} -\frac{2}{\alpha+1}
\sqrt{A+\frac{B\alpha}{\rho_{0} ^{\alpha+1}}} + \frac{\sqrt{A}}{\alpha+1}
\ln\left( 2\sqrt{A}\rho_{0} ^{\alpha+1}
\left( \sqrt{A+\frac{B\alpha}{\rho_{0} ^{\alpha+1}}} + \sqrt{A}\right) +B\alpha \right), ~~\rho_{1} < \rho < \rho_{0}, \\

\lambda_{1}(V_{0}) \leq \lambda_{1}(U) \leq \lambda_{1}(V_{1}),
\end{array}\right.\tag{4.14}
\end{equation}
where $U(t, x) = V(\frac{x}{t})$. Setting $W_{0} =u_{0} -\frac{2}{\alpha+1}
\sqrt{A+\frac{B\alpha}{\rho_{0} ^{\alpha+1}}} + \frac{\sqrt{A}}{\alpha+1}
\ln\left( 2\sqrt{A}\rho_{0} ^{\alpha+1}
\left( \sqrt{A+\frac{B\alpha}{\rho_{0} ^{\alpha+1}}} + \sqrt{A}\right) +B\alpha \right)$ and $N = \sqrt{A+\frac{B\alpha}{\rho ^{\alpha+1}}} $, then we have

\begin{equation}\label{key}
u= \eta+N, ~~~~~~~ \frac{1}{\rho^{\alpha+1}}=\frac{N^{2}-A}{B\alpha}.\tag{4.15}
\end{equation}
  Inserting $u$ into $(4.14)_{2}$, we obtain

\begin{equation}\label{key}
\eta + \frac{(\alpha-1)N}{\alpha+1}+\frac{\sqrt{A}}{\alpha+1}\left[ \ln B\alpha+ \ln \frac{N+\sqrt{A}}{N-\sqrt{A}} \right] = W_{0}. \tag{4.16}
\end{equation}
Further, (4.16) can be rewritten as
\begin{equation}\label{key}
\frac{(\alpha-1)N}{\alpha+1}+\frac{\sqrt{A}}{\alpha+1}\left[  \ln \frac{N+\sqrt{A}}{N-\sqrt{A}} \right] = W_{0}-\eta-\frac{\sqrt{A}\ln B\alpha}{\alpha+1}. \tag{4.17}
\end{equation}
It is too difficult to obtain completely explicit solution $\rho(\eta)$ from (4.17)  due to the structure of implicit equation. Next, we turn to prove the monotonicity of $\rho(\eta)$, which will play an important role in the following discussion. Firstly, we calculate

\begin{equation}\label{key}
N' =\left[  \sqrt{A+\frac{B\alpha}{\rho ^{\alpha+1}}}\right]' = -\frac{(\alpha+1)B\alpha 
	 \rho\prime(\eta)}{2 \sqrt{A+\frac{B\alpha}
		{\rho ^{\alpha+1}}}\rho^{\alpha+2}} = -\frac{(\alpha+1)B\alpha 
	\rho\prime(\eta)}{2N\rho^{\alpha+2}}, \tag{4.18}
\end{equation}
and

\begin{equation}\label{key}
\left[  \ln \frac{N+\sqrt{A}}{N-\sqrt{A}} \right]' = \left[ \ln\left( N+\sqrt{A}\right) - \ln\left( N-\sqrt{A}\right)\right]' = -\frac{2N'\sqrt{A}}{\left( N+\sqrt{A}\right)\left( N-\sqrt{A}\right)}. \tag{4.19}
\end{equation}
Then, substituting  $(4.18)$ into $(4.19)$, we have
\begin{equation}\label{key}
\left[  \ln \frac{N+\sqrt{A}}{N-\sqrt{A}} \right]'  = \frac{(\alpha+1)B\alpha\sqrt{A}\rho\prime(\eta)}{N(N-\sqrt{A})
	(N+\sqrt{A})\rho^{\alpha+2}}.\tag{4.20}
\end{equation}
Next,  performing implicit differentiation with respect to $\eta$ in (4.17) and reorganizing,  we have
\begin{equation}\label{key}
\frac{-B\alpha(\alpha-1)\rho \prime(\eta)}{2N\rho^{\alpha+2}}
+ \frac{A B\alpha  \rho \prime(\eta)}{N(N-\sqrt{A})
	(N+\sqrt{A})\rho^{\alpha+2}} = -1. \tag{4.21}
\end{equation}
 Further simplification based on (4.15) yields
\begin{equation}\label{key}
\frac{-(\alpha-1)\left( N^{2}-A\right) \rho \prime(\eta)}{2N\rho}
+ \frac{ A \rho \prime(\eta)}{N\rho} = -1. \tag{4.22}
\end{equation}
By a simple calculation, (4.22) is equivalent to
\begin{equation}\label{key}
\frac{\left[ (\alpha+1)A+(1-\alpha) N^{2}\right] \rho \prime(\eta)}{2N\rho} = -1. \tag{4.23}
\end{equation}
Eventually, we can get
\begin{equation}\label{key}
	 \rho \prime(\eta) = \frac{-2N \rho(\eta)}{(\alpha+1)A+
	 	(1-\alpha)N^{2}}. \tag{4.24}
\end{equation}
Since $\rho(\eta)>0$, $N= \sqrt{A+\frac{B\alpha}{\rho ^{\alpha+1}}}>0$ and $0<\alpha<1$, we can deduce
\begin{equation}\label{key}
	 \rho \prime(\eta)<0, \tag{4.25}
\end{equation}
which means that $\rho(\eta)$ of the first family rarefaction wave $R_{1}(U(x, t))$
is a monotone decreasing function of $\eta$.

Similarly, 
from $(4.14)_{1}$, we have $\sqrt{A+\frac{B\alpha}{\rho ^{\alpha+1}}} = u-\eta$. Then, substituting it into  $(4.14)_{2}$ immediately yields

\begin{equation}\label{key}
\frac{(\alpha-1)u}{\alpha+1}+\frac{\sqrt{A}}{\alpha+1}\ln
\frac{u-\eta+\sqrt{A}}{u-\eta-\sqrt{A}} = W_{0}-\frac{2\eta}{\alpha+1}-\frac{\sqrt{A}\ln B\alpha}{\alpha+1}. \tag{4.26}
\end{equation}
In view of state $V_{1}=(\rho_{1}, 0)$, $u(\eta)$ satisfies the boundary condition
\begin{equation}\label{key}
u(\eta_{0})=0. \tag{4.27}
\end{equation}
Then, a combination of (4.26) and (4.27) leads to
\begin{equation}\label{key}
\frac{\sqrt{A}}{\alpha+1}\ln
\frac{-\eta_{0}+\sqrt{A}}{-\eta_{0}-\sqrt{A}} = W_{0}-\frac{2\eta_{0}}{\alpha+1}
-\frac{\sqrt{A}\ln B\alpha}{\alpha+1}. \tag{4.28}
\end{equation}

With a simple calculation, (4.28) is equivalent to
\begin{equation}\label{key}
\frac{\eta_{0}-\sqrt{A}}{\eta_{0}+\sqrt{A}}
= \exp\left({ \frac{(\alpha+1)W_{0}-2\eta_{0}-\sqrt{A}\ln B\alpha}
	{\sqrt{A}}}  \right) . \tag{4.29}
\end{equation}
where
\begin{equation}\label{key}
\begin{split}
W_{0}& =u_{0} -\frac{2}{\alpha+1}
\sqrt{A+\frac{B\alpha}{\rho_{0} ^{\alpha+1}}} + \frac{\sqrt{A}}{\alpha+1}
\ln\left( 2\sqrt{A}\rho_{0} ^{\alpha+1}
\left( \sqrt{A+\frac{B\alpha}{\rho_{0} ^{\alpha+1}}} + \sqrt{A}\right) +B\alpha \right) \\
&= -\sqrt{2} - \frac{2\sqrt{A+B \alpha}}{\alpha+1} + \frac{\sqrt{A}}{\alpha+1} \ln \left[2\sqrt{A} \left(\sqrt{A+B \alpha}+\sqrt{A}\right) + B\alpha  \right]\\
& =  -\sqrt{2} - \frac{2\sqrt{2}}{(\alpha+1)M_{0}} +\frac{\sqrt{A}}{\alpha+1}\ln\left(\frac{2\sqrt{2A}M_{0}+AM_{0}^{2}+2}{M_{0}^{2}} \right).
\end{split}. \tag{4.30}
\end{equation}
Taking (4.30) into (4.29), one obtains

\begin{equation}\label{key}
	 \frac{\eta_{0}-\sqrt{A}}{\eta_{0}+\sqrt{A}}
= \exp  \left( \frac{-\sqrt{2}(\alpha+1)-\frac{2\sqrt{2}}{M_{0}}}{\sqrt{A}}\right)\cdot \left( \frac{\sqrt{2}+\sqrt{A}M_{0}}{\sqrt{2}-\sqrt{A}M_{0}}\right)\cdot \exp\left( \frac{-2\eta_{0}}{\sqrt{A}}\right) 
	  . \tag{4.31}
\end{equation}

\noindent Let $C = \exp  \left( \frac{-\sqrt{2}(\alpha+1)-\frac{2\sqrt{2}}{M_{0}}}{\sqrt{A}}\right)\cdot \left( \frac{\sqrt{2}+\sqrt{A}M_{0}}{\sqrt{2}-\sqrt{A}M_{0}}\right)$ and $Q_{0} = \frac{\eta_{0}}{\sqrt{A}}$. It is easy to verify 

\begin{equation}\label{key}
\begin{split}
C& = \exp  \left( \frac{-\sqrt{2}(\alpha+1)-\frac{2\sqrt{2}}{M_{0}}}{\sqrt{A}}\right)\cdot \left( \frac{ \left(\sqrt{2}+\sqrt{A}M_{0}\right)^2}{2-AM_{0}^2
}\right)  \\
&= \exp  \left( \frac{-\sqrt{2}(\alpha+1)-\frac{2\sqrt{2}}{M_{0}}}{\sqrt{A}}\right)\cdot \left( \frac{ \left(\sqrt{2}+\sqrt{A}M_{0}\right)^2}{\alpha BM_{0}^2
}\right)>0.
\end{split}
. \tag{4.32}
\end{equation}
Then, (4.31) can rewritten as
\begin{equation}\label{key}
\exp (-2Q_{0}) = \frac{Q_{0}-1}{C(Q_{0}+1)}. \tag{4.33}
\end{equation}
In what follows, we define 
\begin{equation}\label{key}
f(Q) = \exp (-2Q) - \frac{1}{C} + \frac{2}{C(Q+1)}. \tag{4.34}
\end{equation}

\noindent Then, it follows from (4.33) that there exists a solution $Q_{0} \in (-\infty, 0)$ satisfies $f(Q_{0}) = 0$, which can be verified by  the graph relationships between the exponential function and inverse proportional function. Since (4.33)  is a transcendental equation, it is difficult to find an explicit analytical solution $\eta_{0}$, we turn to study the monotonicity of function (4.34) to demonstrate the existence of a rarefaction wave solution. 
Simple calculations show
\begin{equation}\label{key}
f'(Q) =-2\left( \exp(-2Q)+\frac{1}{C(Q+1)^2}\right)< 0.\tag{4.35}
\end{equation}
\noindent Considering $V_{0} = (1, -\sqrt{2})$, from $(4.14)_1$, 
we have $$\eta_{-\sqrt{2}} = \lambda_{1}(V_{0}) = \sqrt{2}\left( -1 - M_{0}^{-1}\right). $$
\noindent
For $Q_{-\sqrt{2}} = \frac{\eta_{-\sqrt{2}}}{\sqrt{A}} $, it follows from  (4.34) that
\begin{equation}\label{key}
\begin{split}
f(Q_{-\sqrt{2}})&= \exp \left(\frac{2\sqrt{2}\left(1+ \frac{1}{M_{0}} \right) }{\sqrt{A}} \right) - 
 \exp \left(\frac{\sqrt{2}(\alpha+1)+\frac{2\sqrt{2}}{M_{0}}  }{\sqrt{A}} \right) \frac{(\sqrt{2}-\sqrt{A}M_{0})(\sqrt{2}+\sqrt{2}M_{0}
 	+\sqrt{A}M_{0})}
 {(\sqrt{2}+\sqrt{A}M_{0})(\sqrt{2}+\sqrt{2}M_{0}-\sqrt{A}M_{0})}
 \\
&=\exp\left( \frac{\sqrt{2}+\frac{\sqrt{2} }{M_{0}}}{\sqrt{A}}
\right)  \left[ \frac{\exp\left( \frac{\sqrt{2}+\frac{\sqrt{2} }{M_{0}}}{\sqrt{A}}
	\right) (2+2M_{0}-AM_{0}^2+\sqrt{2A}M_{0}^2)
	-\exp\left( \frac{\sqrt{2}\alpha+\frac{\sqrt{2} }{M_{0}}}{\sqrt{A}}
	\right) (2+2M_{0}-AM_{0}^2-\sqrt{2A}M_{0}^2)
}
{(\sqrt{2}+\sqrt{A}M_{0})(\sqrt{2}+\sqrt{2}M_{0}-\sqrt{A}M_{0})}
\right] \end{split} \notag
\end{equation}
Since $\sqrt{A}M_{0}<\sqrt{A+B\alpha}M_{0}=\sqrt{2}$, we have 
\begin{equation}
f(Q_{-\sqrt{2}})>\frac{\exp\left( \frac{\sqrt{2}+\frac{\sqrt{2} }{M_{0}}}{\sqrt{A}}
	\right) 2\sqrt{2A}M_{0}^2}
{(\sqrt{2}+\sqrt{A}M_{0})(\sqrt{2}+\sqrt{2}M_{0}-\sqrt{A}M_{0})}>0,\tag{4.36}
\end{equation}
 i.e., $f(Q_{-\sqrt{2}}) > f(Q_{0})$. From the above discussion, it follows that 
\begin{equation}\label{key}
V_{m}\left( \frac{x}{t} \right) = (\rho(\eta), u(\eta)),\\~~~~~~
\eta_{-\sqrt{2}} \leq \eta \leq \eta_{0},\\~~~~~~
0 < \rho(\eta_{0}) \leq \rho(\eta) \leq \rho(\eta_{-\sqrt{2}}). \tag{4.37}
\end{equation}
Therefore,  the first family rarefaction wave $R_{1}(U)$ is a solution of (2.8), (2.6) and (2.14) connecting $V_{0} = (1, -\sqrt{2})$, $V_{1} = (\rho_{1}, 0)$ in the domain $\Omega$.

Analogously, for the second family rarefaction wave $R_{2}\left( U(x, t)\right) $, we can deduce  the self-similar solution  as follow: 

\begin{equation}\label{key}
	\left\{\begin{array}{ll} 
		\xi = \frac{x}{t} = \lambda_{2}(U) = u + \sqrt{A + \frac{B \alpha}{\rho^{\alpha+1}}},\\
		
		u +\frac{2}{\alpha+1}
		\sqrt{A+\frac{B\alpha}{\rho^{\alpha+1}}} - \frac{\sqrt{A}}{\alpha+1}
		\ln\left( 2\sqrt{A}\rho^{\alpha+1}
		\left( \sqrt{A+\frac{B\alpha}{\rho ^{\alpha+1}}} + \sqrt{A}\right) +B\alpha \right)  =\\
		u_{0} +\frac{2}{\alpha+1}
		\sqrt{A+\frac{B\alpha}{\rho_{0} ^{\alpha+1}}} - \frac{\sqrt{A}}{\alpha+1}
		\ln\left( 2\sqrt{A}\rho_{0} ^{\alpha+1}
		\left( \sqrt{A+\frac{B\alpha}{\rho_{0} ^{\alpha+1}}} + \sqrt{A}\right) +B\alpha \right), ~~\rho_{1} < \rho < \rho_{0}, \\
		
		\lambda_{2}(V_{0}) \leq \lambda_{2}(U) \leq \lambda_{2}(V_{1}).
	\end{array}\right.\tag{4.38}
\end{equation}
By the same analysis of the first family rarefaction wave, we can get $$\xi_{-\sqrt{2}} = \lambda_{2}(V_{0}) = \sqrt{2}\left( -1 + M_{0}^{-1}\right),$$
 and
\begin{equation}\label{key}
	\frac{(\alpha-1)N}{\alpha+1}+\frac{\sqrt{A}}{\alpha+1}\left[  \ln \frac{N+\sqrt{A}}{N-\sqrt{A}} \right] = \xi-W_{1}-\frac{\sqrt{A}\ln B\alpha}{\alpha+1}. \tag{4.39}
\end{equation}
where 

\begin{equation}\label{key}
\begin{split}
W_{1} &= 	u_{0} +\frac{2}{\alpha+1}
\sqrt{A+\frac{B\alpha}{\rho_{0} ^{\alpha+1}}} - \frac{\sqrt{A}}{\alpha+1}
\ln\left( 2\sqrt{A}\rho_{0} ^{\alpha+1}
\left( \sqrt{A+\frac{B\alpha}{\rho_{0} ^{\alpha+1}}} + \sqrt{A}\right) +B\alpha \right)\\
&=
-\sqrt{2} + \frac{2\sqrt{2}}{\left(\alpha+1  \right) M_{0}} - \frac{\sqrt{A}}{\left(\alpha+1  \right) }\ln\left(\frac{2\sqrt{2A}M_{0}+AM_{0}^{2}+2}{M_{0}^{2}} \right).
\end{split}
\tag{4.40}
\end{equation}

Differentiating the both sides of (4.39) with respect to $\xi$, we obtain
\begin{equation}\label{key}
-\frac{(\alpha-1)\left(N^{2}-A \right)\rho'(\xi) }{2N\rho}
+\frac{A\rho'(\xi)}{N\rho}
=1. \tag{4.41}
\end{equation}
Then, rearranging (4.41) leads to
\begin{equation}\label{key}
\rho'(\xi) = \frac{2N\rho}{2A-\left(\alpha-1 \right)\left(N^2-A \right)  }
= \frac{2N\rho}{A+\left(1-\alpha\right)N^2+\alpha A}>0, \tag{4.42}
\end{equation}
 which means that the density of the second family rarefaction wave $R_{2}(U(x, t))$ is monotonic increasing with respect to $\xi$.
 
Based on $(4.38)_1$, taking $u=\xi-\sqrt{A + \frac{B \alpha}{\rho^{\alpha+1}}}$ into $(4.38)_2$ to deduce

\begin{equation}\label{key}
\frac{(\alpha-1)u}{\alpha+1}-\frac{\sqrt{A}}{\alpha+1}\ln
\frac{u-\xi-\sqrt{A}}{u-\xi+\sqrt{A}} = W_{1}-\frac{2\xi}{\alpha+1}
+\frac{\sqrt{A}\ln B\alpha}{\alpha+1}. \tag{4.43}
\end{equation}
In view of state $V_{1} = (\rho_{1}, 0)$, $u(\xi) $ satisfies the boundary condition $u(\xi_{0}) = 0$. Inserting $u(\xi_{0})$ into (4.43) to obtain 

\begin{equation}\label{key}
	\frac{\xi_{0}+\sqrt{A}}{\xi_{0}-\sqrt{A}}
	= \exp\left({ \frac{2\xi_{0}-(\alpha+1)W_{1}-\sqrt{A}\ln B\alpha}
		{\sqrt{A}}}  \right) . \tag{4.44}
\end{equation}
It follows from (4.44) and (4.40) that

 \begin{equation}\label{key}
 	\frac{\xi_{0}+\sqrt{A}}{\xi_{0}-\sqrt{A}}
 	= \exp  \left( \frac{\sqrt{2}(\alpha+1)-\frac{2\sqrt{2}}{M_{0}}}{\sqrt{A}}\right)\cdot \left( \frac{2\sqrt{2A}M_{0}+AM_{0}^2+2}{\alpha BM_{0}^2}\right)\cdot \exp\left( \frac{2\xi_{0}}{\sqrt{A}}\right) 
 	. \tag{4.45}
 \end{equation}
Let $C_{1} = \exp  \left( \frac{\sqrt{2}(\alpha+1)-\frac{2\sqrt{2}}{M_{0}}}{\sqrt{A}}\right)\cdot \left( \frac{2\sqrt{2A}M_{0}+AM_{0}^2+2}{\alpha BM_{0}^2}\right)$ and $Q=\frac{\xi}{\sqrt{A}}$. Then, we  define $g(Q) = \exp (2Q) - \frac{1}{C_{1}} - \frac{2}{C_{1}(Q - 1)}$ with $g(Q_{0})=g\left( \frac{\xi_{0}}{\sqrt{A}}\right)  = 0$. Simple calculations lead to 
\begin{equation}\label{key}
g'(Q) = 2\exp(2Q)+\frac{2}{C_{1}\left(Q-1 \right)^2 }> 0,\tag{4.46}
\end{equation}
and
\begin{equation}\label{key}
\begin{split}
g(Q_{-\sqrt{2}})&=g\left(\frac{\xi_{-\sqrt{2}}}{\sqrt{A}} \right)=\exp\left( \frac{-\sqrt{2}+\frac{\sqrt{2} }{M_{0}}}{\sqrt{A}} \right) \\
&\times
 \left[\exp\left( \frac{-\sqrt{2}+\frac{\sqrt{2} }{M_{0}}}{\sqrt{A}}
 \right)- \frac{ 
	\exp\left( \frac{-\sqrt{2}\alpha+\frac{\sqrt{2} }{M_{0}}}{\sqrt{A}}
	\right)\left( \frac{\sqrt{2}}{M_{0}}-\sqrt{2}+\sqrt{A}\right) (\sqrt{2}-\sqrt{A}M_{0})
}
{\left( \frac{\sqrt{2}}{M_{0}}-\sqrt{2}-\sqrt{A}\right) (\sqrt{2}+\sqrt{A}M_{0})}
\right] 
\end{split}\tag{4.47}
\end{equation}
If $\frac{\sqrt{2}}{M_{0}}-\sqrt{2}-\sqrt{A} >0$, we have 
\begin{equation}\label{key}
g(Q_{-\sqrt{2}})<\exp\left( \frac{-\sqrt{2}+\frac{\sqrt{2} }{M_{0}}}{\sqrt{A}}\right) 
\left[ \frac{\exp\left( \frac{-\sqrt{2}+\frac{\sqrt{2} }{M_{0}}}{\sqrt{A}}\right) (-2\sqrt{2A}M_{0})}
{\left( \frac{\sqrt{2}}{M_{0}}-\sqrt{2}-\sqrt{A}\right) (\sqrt{2}+\sqrt{A}M_{0})}
\right] 
<0.
\tag{4.48}
\end{equation}
 which  infers $\xi_{-\sqrt{2}} < \xi_{0}$ and $\rho(\xi_{-\sqrt{2}}) < \rho(\xi_{0})$,  contradicting with the requirement $\rho_{0}>\rho_{1}$. 
 
\noindent Similarly, if  $\frac{\sqrt{2}}{M_{0}}-\sqrt{2}-\sqrt{A} <0$, we can obtain
 \begin{equation}\label{key}
 \begin{split}
 g(Q_{-\sqrt{2}})	>0
 \end{split}. \tag{4.49}
 \end{equation}
 which means $\xi_{-\sqrt{2}} > \xi_{0}$,  contradicting with the requirement
 	$\lambda_{2}(V_{0}) \leq \lambda_{2}(U) \leq \lambda_{2}(V_{1})$.  Hence, the second family rarefaction wave $R_{2}(U(x, t))$ is not a physical solution as the piston recedes from the gas. 

Consequently, we have the following result.

\noindent \textbf{Theorem 4.2.} \label{{Theorem 4.2}}~~ For $v_{0} = \sqrt{2}$ and any fixed $A$ and $B$, as the piston recedes from the MCG for any $M_{0} \in (0, \infty)$,  only one physically relevant rarefaction wave solution will be present in front of the piston to connect two states $V_{0} = (1, -\sqrt{2})$ and $V_{1} = (\rho_{1}, 0)$ in the domain $\Omega$, and the concentration will never occur. 

\noindent \textbf{Remark 4.2.}~~  The results obtained here, to some extent, not only extend the results and proofs in [10, 17], in which  the pressure state only contains one parameter, but also cover the results obtained in [10] for generalized Chaplygin gas as $A \rightarrow 0$, which will be analyzed and identified in Section 5.

\section{Limit behavior of the piston solutions as A $\rightarrow$ 0}
~~~~ This section discusses the limit behavior of piston solutions to (2.8) and (1.5) as the pressure (1.5) approaches the generalized Chaplygin gas (1.3), i.e., $A \rightarrow 0$. 
Through the analysis in Section 2, taking the limit  $A \rightarrow 0$, the parameter B in (2.10) approaches to 
\begin{equation}\label{key}
B = \frac{2}{\alpha M_{0}^{2}}, \tag{5.1}
\end{equation}
and the pressure in (2.11) tends to 
\begin{equation}\label{key}
P(\rho) = -\frac{2}{\alpha M_{0}^{2}\rho^{\alpha }}. \tag{5.2}
\end{equation}
Both (5.1) and (5.2) play an important role in forming singularity as the piston moves forward to the MCG and $A\rightarrow 0$. From (5.1) and (5.2), similar with the proof of the theorem 4.1, we can obtain
\begin{equation}\label{key}
\lim\limits_{A \rightarrow 0} M_{0}^{2} = {\alpha}^{-1} \left(1-\frac{1}{\rho_{1}} \right) \left(1-\frac{1}{\rho_{1}^{\alpha}} \right) . \tag{5.3}
\end{equation}
From (5.3), it is easy to discover, if $0<M_{0}< \sqrt{\alpha^{-1}}$, it follows that $\rho_{1}>1$, which shows that there exists an integral weak solution in front of the piston. While $M_{0}\geq \sqrt{\alpha^{-1}}$, it contradicts with $\rho_{1}>1$, i.e., the piston solution cannot be constructed by using shock waves.

Then, we can construct a Radon measure solution based on definition 3.2 for the case $M_{0}\geq \sqrt{\alpha^{-1}}$,
\begin{equation}\label{key}
\tau =  I_{\Omega}L^{2} + w_{\rho}(t)\delta_{\{x=0,~ t \geq 0\} }, 
~~~m =  \sqrt{2} I_{\Omega}L^{2} , ~~~n = 2 I_{\Omega}L^{2} ,  ~~~\wp = -\frac{2}{\alpha M_{0}^2}I_{\Omega}L^{2},\tag{5.4}
\end{equation}
where $I_{\Omega}$ is the characteristic function of $\Omega$. A Radon measure solution is a generalization of an ordinary shock wave. Informally speaking, it is a kind of discontinuity, on which at least one of the state variables may develop an extreme concentration in the form of a weighted Dirac delta function with the discontinuity as its support. It is more compressive than an ordinary shock wave in the sense of distributions.

Substituting (5.4) into (3.5), we infer that 
\begin{equation}\label{key}
\begin{split}
0 &= <\tau, \partial t \phi> + <m, \partial x \phi> + 
\int_{-\infty}^{0} \rho_{0} \phi(0, x)dx \\&= \int_{\Omega} \partial t \phi dxdt + \int_{0}^{\infty} w_{\rho}(t) \partial t  \phi(t, 0)dt
+ \sqrt{2}  \int_{0}^{\infty}  \int_{-\infty}^{0}  \partial x  \phi dxdt 
+ \int_{-\infty}^{0} \phi(0, x)dx 
\\&= \int_{-\infty}^{0} \phi(t, x)|_{t=0}^{\infty} dx
+ w_{\rho}(t) \phi(t, 0)|_{t=0}^{\infty} 
-  \int_{0}^{\infty} w'_{\rho}(t) \phi(t, 0)dt + 
\sqrt{2}  \int_{0}^{\infty} \phi(t, x)|_{x=-\infty}^{0} dt +  \int_{-\infty}^{0} \phi(0, x)dx 
\\&=  \int_{0}^{\infty} (\sqrt{2} - w'_{\rho}(t))\phi(t, 0)dt
-w_{\rho}(0)\phi(0, 0).
\end{split} \tag{5.5}
\end{equation}
Due to the arbitrariness of $\phi$, we have
\begin{equation}\label{key}
\left\{\begin{array}{ll} 
w'_{\rho}(t) = \sqrt{2}, ~~~t\geq0,\\
w_{\rho}(0) = 0.
\end{array}\right. \tag{5.6}
\end{equation}
It follows that  
\begin{equation}\label{key}
w_{\rho}(t) = \sqrt{2}t. \tag{5.7}
\end{equation}
Similarly, substituting (5.4) into  (3.6), we have
\begin{equation}\label{key}
\begin{split}
0 &= <m, \partial t \phi> + <n, \partial x \phi> + 
<\wp, \partial x \phi> - <w_{p}(t)\delta_{\{x=0,~ t \geq 0\} }, \phi>
+ \int_{-\infty}^{0} \rho_{0}u_{0} \phi(0, x)dx 
\\&= \sqrt{2}  \int_{0}^{\infty}  \int_{-\infty}^{0}
\partial t \phi dxdt +
2 \int_{0}^{\infty}  \int_{-\infty}^{0}
\partial x \phi dxdt - 
\frac{2}{\alpha M_{0}^2} \int_{0}^{\infty}  \int_{-\infty}^{0}
\partial x \phi dxdt-
 \int_{0}^{\infty} w_{p}(t)  \phi(t, 0)dt
+ \sqrt{2}   \int_{-\infty}^{0} \phi(0, x)dx 
\\&= \int_{0}^{\infty} \left[ 2-w_{p}(t) - \frac{2}{\alpha M_{0}^2}\right] \phi(t, 0)dt.
\end{split} \tag{5.8}
\end{equation}
Since the arbitrariness of $\phi$ and $M_{0}\geq \sqrt{\alpha^{-1}}$, then
\begin{equation}\label{key}
w_{p}(t) = 2 - \frac{2}{\alpha M_{0}^2} \geq 0. \tag{5.9}
\end{equation}
In summary, 
we  conclude that, as $A \rightarrow 0$ and $M_{0} \geq  \sqrt{\alpha^{-1}}$, the shock wave and all gas between the shock wave and  the piston adhere to the piston and then form a concentration of mass like a Dirac measure,  which indicates that 
 the radon measure solution (5.4) is a reasonable solution for the limiting piston problem of modified  Chaplygin gas as $A \rightarrow 0$ and $M_{0} \geq  \sqrt{\alpha^{-1}}$.

For  the receding case, as $A \rightarrow 0$,  (4.10), (4.11) have changed to

\begin{equation}\label{key}
\lambda_{1}(U) = u - \sqrt{\frac{B \alpha}{\rho^{\alpha+1}}}, ~~~~~ \lambda_{2}(U) = u +  \sqrt{\frac{B \alpha}{\rho^{\alpha+1}}}, \tag{5.10}
\end{equation}
with the corresponding right eigenvectors

\begin{equation}\label{key}
\vec{r}_{1} = \left(  \sqrt{\frac{B \alpha}{\rho^{\alpha+1}}}, -\rho\right)^{T},~~~~~
\vec{r}_{2} = \left(  \sqrt{ \frac{B \alpha}{\rho^{\alpha+1}}}, \rho\right)^{T}, \tag{5.11}
\end{equation} 
satisfying 
\begin{equation}\label{key}
\triangledown \lambda_{i}(U) \cdot\vec{r}_{i}  \neq 0, ~~~ ~(i=1, 2,~~~~ 0<\alpha<1). \tag{5.12}
\end{equation}
Then (4.14) has degenerated into 
\begin{equation}\label{key}
\left\{\begin{array}{ll} 
\eta = \frac{x}{t} = \lambda_{1}(U) = u -   \sqrt{\frac{B \alpha}{\rho^{\alpha+1}}},\\

u -\frac{2}{\alpha+1}
\sqrt{\frac{B\alpha}{\rho^{\alpha+1}}}   =
u_{0} -\frac{2}{\alpha+1}
\sqrt{\frac{B\alpha}{\rho_{0} ^{\alpha+1}}} , ~~\rho_{1} < \rho < \rho_{0}, \\

\lambda_{1}(V_{0}) \leq \lambda_{1}(U) \leq \lambda_{1}(V_{1}). 
\end{array}\right.\tag{5.13}
\end{equation}
We can observe that  the self-similar solution above is totally the same as the Eq. (3.8) in [10]. So, performing the similar analysis in [10], we also  obtain only the first rarefaction wave $R_{1}\left( U(x,t)\right) $ is a physical solution. Readers can refer to [10] for more details.

In particular, when $\alpha=1$, we get $ \triangledown \lambda_{i}(U) \cdot\vec{r}_{i} = 0 ~ (i=1, 2)$, which means that both the characteristic fields are linearly degenerate and the elementary waves only involve contact discontinuities, i.e.,  the rarefaction waves coincide with that of the shock waves in the physical plane. Replacing $\rho_{0}, u_{0} $ by  $\rho_{0} = 1, u_{0} = -\sqrt{2}$  and the $(4.3)_{1}$,
 we have
\begin{equation}\label{key}
\sigma = \frac{\sqrt{2}}{\rho_{1}-1}. \tag{5.14}
\end{equation}
Then, by $(4.3)_{2}$,  (5.2) and (5.14), we have
\begin{equation}\label{key}
(2+P_{0})\rho_{1}^{2}-2P_{0}\rho_{1}+P_{0} = 0. \tag{5.15} 
\end{equation}
For $M_{0} \in (0, \infty)$ and the non-negativeness of $\rho_{1}$, it follows from (5.15) that

\begin{equation}\label{key}
 \rho_{1} = \frac{2P_{0}+2\sqrt{2}\sqrt{-P_{0}}}{4+2P_{0}}
 =\frac{1}{1+M_{0}}<1, \tag{5.16} 
\end{equation}
i.e., we prove the existence of rarefaction wave solutions when the piston recedes from the MCG  and $A \rightarrow 0$.

Therefore, 
as $A \rightarrow 0$ and $0<\alpha \leq 1$, the limit of solution of 
 piston problem  (2.8) and (1.5)  is 
similar with that of  generalized or pure Chaplygin Euler equations [10, 17].

 \sec{\Large\bf    Acknowledgments}
 This work is supported by the  Minnan Normal University (Grant No. KJ2021020)  and the Department of Education, Fujian Province (Grant No. JAT210254).

 \sec{\Large\bf    Author Declarations}
  \subsec{\large\bf    Conflict of Interest}
  The authors have no conflicts to disclose.
  
   \sec{\Large\bf    Data Availability}
Data sharing is not applicable to this article as no new data were created or analyzed in this study.

\bibliography{aipsamp}

\end{document}